\documentclass[11pt]{amsart}
\usepackage{amssymb}

\newtheorem{thm}{Theorem}[section]

\newtheorem{lem}[thm]{Lemma}

\theoremstyle{definition}

\newtheorem{rem}[thm]{Remark}

\numberwithin{equation}{section}

\begin{document}
\title[The number of integer points]{The number of integer points in a family of anisotropically expanding domains}

\author[Y. A. Kordyukov]{Yuri A. Kordyukov}
\address{Institute of Mathematics\\
         Russian Academy of Sciences\\
         112~Chernyshevsky str.\\ 450008 Ufa\\ Russia} \email{yurikor@matem.anrb.ru}

\author[A. A. Yakovlev]{Andrey A. Yakovlev}
\address{Institute of Mathematics\\
         Russian Academy of Sciences\\
         112~Chernyshevsky str.\\ 450008 Ufa\\ Russia} \email{yakovlevandrey@yandex.ru}

\date{}

\begin{abstract}
We investigate the remainder in the asymptotic formula for the
number of integer points in a family of bounded domains in the
Euclidean space, which remain unchanged along some linear subspace
and expand in the directions, orthogonal to this subspace. We prove
some estimates for the remainder, imposing additional assumptions on
the boundary of the domain. We study the average remainder
estimates, where the averages are taken over rotated images of the
domain by a subgroup of the group ${\rm SO}(n)$ of orthogonal
transformations of the Euclidean space ${\mathbb R}^n$.

Using these results, we improve the remainder estimate in the
adiabatic limit formula for the eigenvalue distribution function of
the Laplace operator associated with a bundle-like metric on a
compact manifold equipped with a Riemannian foliation in the
particular case when the foliation is a linear foliation on the
torus and the metric is the standard Euclidean metric on the torus.
\end{abstract}
\maketitle

\section{Preliminaries and main results}
\subsection{The setting of the problem}
A classical problem on integer points distribution consists in the
study of the asymptotic behavior of the number of points of the
integer lattice ${\mathbb Z}^n$ in a family of homothetic domains in
${\mathbb R}^n$. This problem is originated in the Gauss problem on
the number of integer points in the disk, where it is directly
related with the arithmetic problem on the number of representations
of an integer as a sum of two squares, and sufficiently well studied
(see, for instance, books
\cite{Fricker,Gruber-Lekker,Huxley,Kraetzel} and the references
therein).

In this paper we investigate much less studied problem of counting
integer points in a family of anisotropically expanding domains.
More precisely, let $F$ be a $p$-dimensional linear subspace of
$\mathbb {R}^n$ and $H=F^{\bot}$ the $q$-dimensional orthogonal
complement of $F$ with respect to the standard inner product
$(\cdot,\cdot)$ in $\mathbb {R}^n$, $p+q=n$. For any
$\varepsilon>0$, consider the linear transformation $T_\varepsilon :
\mathbb {R}^n\to \mathbb {R}^n$ given by
\[
T_\varepsilon(x)=\begin{cases} x, & \text{if}\ x\in F, \\
\varepsilon^{-1}x, & \text{if}\ x\in H.
\end{cases}
\]
For any bounded set $S$ in $\mathbb{R}^n$, we put
\begin{equation}\label{e:nS}
n_\varepsilon (S)=\# (T_\varepsilon(S)\cap \mathbb{Z}^n), \quad
\varepsilon>0.
\end{equation}

The study of the asymptotic behavior of $n_\varepsilon (S)$ as
$\varepsilon \to 0$ for general domains in ${\mathbb R}^n$ was
started in \cite{lattice-points}. In particular, the following
asymptotic formula has been proved.

Let $\Gamma = \mathbb{Z}^n\cap F$. $\Gamma$ is a free abelian group.
Denote by $r=\operatorname{rank} \Gamma \leq p$ the rank of
$\Gamma$. Let $V$ be the $r$-dimensional subspace of $\mathbb R^n$
spanned by the elements of $\Gamma$. Observe that $\Gamma$ is a
lattice in $V$. Let $\Gamma^*$ denote the lattice in $V$, dual to
the lattice $\Gamma$:
\[
\Gamma^*=\{ \gamma^*\in V :(\gamma^*,\Gamma)\subset\mathbb Z\}.
\]
For any $x\in V$, we denote by $P_{x}$ the $(n-r)$-dimensional
affine subspace of $\mathbb {R}^n$, passing through $x$ orthogonal
to $V$.

\begin{thm}[\cite{lattice-points}, Theorem 1.1]\label{t:lattice_points0} For any bounded open
set $S$ in $\mathbb{R}^n$ with smooth boundary, we have
\begin{equation}\label{f:lattice_points1}
 n_\varepsilon (S) =
\frac{\varepsilon^{-q}}{{\rm vol}(V/\Gamma)}
 \sum_{\gamma^*\in\Gamma^*} {\rm vol}_{n-r} (P_{\gamma^*}\cap S)
+O(\varepsilon^{\frac{1}{p-r+1}-q}), \quad \varepsilon \to 0.
\end{equation}
\end{thm}
Here, if $\Gamma$ is trivial, by definition we put ${\rm
vol}(V/\Gamma)=1$.

Actually, as it can be easily seen from the proof of this theorem,
it holds for any set $S\subset \mathbb{R}^n$ such that for any
$\gamma^*\in\Gamma^*$ the intersection $P_{\gamma^*} \cap S$ is a
bounded open set in $P_{\gamma^*}$ with Lipschitz boundary. Using
this fact, one can easily prove a more general form of the formula
\eqref{f:lattice_points1}.

\begin{thm}\label{t:lattice_points-basic} For any
set $S$ in $\mathbb{R}^n$ such that for any $\gamma^*\in\Gamma^*$
the intersection $P_{\gamma^*} \cap S$ is a bounded open set in
$P_{\gamma^*}$ which is Jordan measurable, we have
\begin{equation}\label{f:lattice_points-basic}
 n_\varepsilon (S) \sim
\frac{\varepsilon^{-q}}{{\rm vol}(V/\Gamma)}
 \sum_{\gamma^*\in\Gamma^*} {\rm vol}_{n-r} (P_{\gamma^*}\cap S), \quad \varepsilon \to 0.
\end{equation}
\end{thm}

In this paper, we continue the study of the remainder in the formula
\eqref{f:lattice_points-basic} given by
\[
R_\varepsilon (S)=n_\varepsilon (S) - \frac{\varepsilon^{-q}}{{\rm
vol}(V/\Gamma)}
 \sum_{\gamma^*\in\Gamma^*} {\rm vol}_{n-r} (P_{\gamma^*}\cap S).
\]

Remark that, in a slightly different context, the problem of
counting integer points in a family of anisotropically expanding
domains was also studied in considerable detail in
\cite{Skriganov89,Skriganov94,Nikichine-Skriganov95,Nikichine-Skriganov98}
(see also the introduction of \cite{lattice-points}).

\subsection{Main results}
The first goal of the paper is to obtain more precise estimates for
$R_\varepsilon (S)$, imposing some additional assumptions on $S$.

\begin{thm}\label{t:lattice_points}
Let $S$ be any subset of $\mathbb{R}^n$ such that for any
$\gamma^*\in\Gamma^*$, the intersection $P_{\gamma^*} \cap S$ is a
bounded open set in $P_{\gamma^*}$ with smooth boundary.

(1) If, for any $\gamma^*\in\Gamma^*$ and $x\in F\cap V^\bot$, the
intersection $S\cap \{\gamma^*+x+H\}$ is strictly convex, then we
have
\begin{equation}\label{f:lattice_points0}
R_\varepsilon (S)=O(\varepsilon^{\frac{2q}{q+1+2(p-r)}-q}), \quad
\varepsilon \to 0.
\end{equation}

(2) If, for any $\gamma^*\in\Gamma^*$, the intersection
$P_{\gamma^*} \cap S$ is strictly convex, then:
\begin{equation}\label{f:lattice_points}
R_\varepsilon (S)= O(\varepsilon^{\frac{2q}{n-r+1}-q}), \quad
\varepsilon \to 0.
\end{equation}
\end{thm}

Remark that the first statement of Theorem~\ref{t:lattice_points} is
a slight improvement of Theorem 1.2 in \cite{lattice-points}.

In \cite{Randol69a,Randol69b}, Randol suggested in the case of a
family of homothetic domains in ${\mathbb R}^n$ to consider instead
of the remainder for a domain $S$ its averages over rotated or over
rotated and translated images of $S$. He observed that estimates for
averages can be substantially smaller than individual estimates. In
\cite{Randol69a,Randol69b}, such results were proved for bounded
convex open sets with analytic boundary. The results of Randol were
extended by Varchenko \cite{Varchenko83} to arbitrary bounded open
sets with smooth boundary, proving a conjecture by Arnold. Finally,
in \cite{Brandolini} the average remainder estimates were proved for
any bounded convex open set $S$ or in the case when the boundary of
$S$ is $C^{3/2}$. We refer to \cite{Brandolini,Gruber-Lekker} for
more information and references on this problem.

The second goal of the paper is to obtain similar results in the
case under consideration. So we study averages of the remainder $
R_\varepsilon (S)$ taken over rotated images of $S$ by a group of
orthogonal transformations of the Euclidean space ${\mathbb R}^n$.
We will consider several subgroups of ${\rm SO}(n)$.

First, we consider the group ${\rm SO}(H)$ of orthogonal
transformations of ${\mathbb R}^n$, which fix any vector of $F$
(and, as a consequence, take $H$ to itself):
\[
{\rm SO}(H)=\{A\in {\rm SO}(n) : A\left|_F\right.={\rm Id}\}.
\]
We will denote by $dh$ the Haar measure on ${\rm SO}(n)$ (and on any
subgroup of ${\rm SO}(n)$ as well).

\begin{thm}\label{t:VarchenkoH}
For any subset $S$ of $\mathbb{R}^n$ such that for any
$\gamma^*\in\Gamma^*$ and $x\in F\cap V^\bot$, the intersection
$S\cap \{\gamma^*+x+H\}$ is a bounded open set in $P_{\gamma^*}$
such that either $S\cap \{\gamma^*+x+H\}$ is convex or the boundary
of $S\cap \{\gamma^*+x+H\}$ is $C^{3/2}$, we have
\[
\int_{SO(H)} |R_\varepsilon
(hS)|\,dh=O(\varepsilon^{\frac{2q}{q+1+2(p-r)}-q}), \quad
\varepsilon \to 0.
\]
\end{thm}

Next, we consider the group ${\rm SO}(V^\bot)$ of orthogonal
transformations of ${\mathbb R}^n$, which fix any vector of $V$:
\[
{\rm SO}(V^\bot)=\{A\in {\rm SO}(n) : A\left|_V\right.={\rm Id}\}.
\]

\begin{thm}\label{t:Varchenko-full}
For any subset $S$ of $\mathbb{R}^n$ such that for any
$\gamma^*\in\Gamma^*$, the intersection $P_{\gamma^*} \cap S$ is a
bounded open set in $P_{\gamma^*}$ with smooth boundary, we have
\[
\int_{SO(V^\bot)} |R_\varepsilon (hS)|\,dh
=O(\varepsilon^{\frac{2q}{n-r+1}-q}), \quad \varepsilon \to 0.
\]
\end{thm}

Finally, Theorem~\ref{t:Varchenko-full} easily implies the result on
the average remainder estimates for the full group ${\rm SO}(n)$.

\begin{thm}\label{t:Varchenko-full2}
For any bounded open set $S$ in $\mathbb{R}^n$ with smooth boundary,
we have
\[
\int_{SO(n)} |R_\varepsilon (hS)|\,dh
=O(\varepsilon^{\frac{2q}{n-r+1}-q}), \quad \varepsilon \to 0.
\]
\end{thm}

\begin{rem}
It appears that, using the results of \cite{Brandolini}, one can
relax the assumptions on smoothness of the boundary of $S$ in
Theorems~\ref{t:Varchenko-full} and \ref{t:Varchenko-full2}.
\end{rem}

\subsection{Applications to adiabatic limits} It is well known that
the Gauss problem on counting integer points in the disk is
equivalent to the problem on the asymptotic behavior of the
eigenvalue distribution function of some elliptic differential
operator on a compact manifold, namely, of the Laplace operator on a
torus. In the case under consideration, there is also an equivalent
asymptotic spectral problem, namely, the problem on the asymptotic
behavior of the eigenvalue distribution function of the Laplace
operator on a torus in the adiabatic limit associated with a linear
foliation.

As above, let $F$ be a $p$-dimensional linear subspace of $\mathbb
{R}^n$ and $H=F^{\bot}$. Consider the $n$-dimensional torus $\mathbb
T^n=\mathbb {R}^n/\mathbb{Z}^n$. Let $\mathcal F$ be the associated
linear foliation on $\mathbb T^n$: the leaf $L_x$ of $\mathcal F$
through $x\in \mathbb T^n$ has the form:
\[
L_x=x+F \mod \mathbb{Z}^n.
\]
The decomposition of $\mathbb{R}^n$ into the direct sum of subspaces
$\mathbb{R}^n=F\oplus H$ induces the decomposition $g=g_{F}+g_{H}$
of the standard Euclidean metric $g$ on $\mathbb{R}^n$ into the sum
of the tangential and transversal components. Define a one-parameter
family $g_{\varepsilon}$ of Euclidean metrics on $\mathbb{R}^n$ by
\begin{equation*}
g_{\varepsilon}=g_{F} + {\varepsilon}^{-2}g_{H}, \quad \varepsilon >
0.
\end{equation*}
We will also consider the metrics $g_{\varepsilon}$ as Riemannian
metrics on $\mathbb{T}^n$.

Let $A=(a_1,\ldots,a_n)\in \mathbb{R}^n$. For any $\varepsilon>0$,
consider the operator $H_\varepsilon$ in $C^\infty(\mathbb T^n)$
defined in the standard linear coordinates $(x_1,x_2,\ldots,x_n)$ in
$\mathbb{R}^n$ by
\[
H_\varepsilon =\sum_{j,\ell=1}^n g_\varepsilon^{j\ell}
\left(\frac{\partial}{\partial x_j} -2\pi
ia_j\right)\left(\frac{\partial}{\partial x_\ell} - 2\pi
ia_\ell\right),
\]
where $g_\varepsilon^{j\ell}$ are the elements of the inverse matrix
of $g_\varepsilon $. The operator $H_\varepsilon$ can be considered
as the magnetic Schr\"odinger operator on the torus $\mathbb T^n$,
associated with the metric $g_\varepsilon$ and the constant magnetic
potential ${\mathbf A}=\sum_{j=1}^n a_jdx_j$. It has a complete
orthogonal systems of eigenfunctions
\[
U_{k}(x)=e^{2\pi i (k,x)}, \quad x\in \mathbb R^n,\quad k\in \mathbb
Z^n,
\]
with the corresponding eigenvalues
\[
\lambda_{k}=(2\pi)^2
\|k-A\|^2_{g_\varepsilon^{-1}}=(2\pi)^2\sum_{j,\ell=1}^n
g_\varepsilon^{j\ell} \left(k_j -a_j\right) \left(k_\ell -
a_\ell\right).
\]

Denote by $N_\varepsilon(\lambda)$ the eigenvalue distribution
function of $H_\varepsilon$:
\[
N_\varepsilon(\lambda)=\sharp\{k\in \mathbb Z^n:\lambda_{k}<\lambda
\}, \quad \lambda \in {\mathbb R}.
\]
It is easy to see that
\[
n_\varepsilon (B_{\sqrt{\lambda}}(A))=N_\varepsilon (4\pi^2\lambda),
\quad \lambda \in {\mathbb R}.
\]
Thus, the problem on the asymptotic behavior of the number
$n_\varepsilon (B_{\sqrt{\lambda}}(A))$ of integer points in the
ellipsoid $T_\varepsilon(B_{\sqrt{\lambda}}(A))$ as $\varepsilon\to
0$ is equivalent to the problem on the asymptotic behavior of the
eigenvalue distribution function $N_\varepsilon(\lambda)$ as
$\varepsilon\to 0$. The limiting procedure $\varepsilon\to 0$ is
often called passing to adiabatic limit. This notion was introduced
by Witten in 1985 in the study of global anomalies in string theory.
We refer the reader to a survey paper \cite{bedlewo2-andrey} for
some historic remarks and references.

In~\cite{adiab} (see also \cite{asymp}), the first author computed
the leading term of the asymptotics of the eigenvalue distribution
function of the Laplace operator associated with a bundle-like
metric on a compact manifold equipped with a Riemannian foliation,
in adiabatic limit. The linear foliation on the torus is a
Riemannian foliation, and a Euclidean metric on the torus is
bundle-like.

A more precise estimate of the remainder in the asymptotic formula
of~\cite{adiab} in this particular case was obtained in
\cite{lattice-points}. As a straightforward consequence of
Theorem~\ref{t:lattice_points}, we improve the remainder estimate
of~\cite{lattice-points}.

\begin{thm}\label{t:lattice_points2}
For $\lambda>0$, the following asymptotic formula holds as
$\varepsilon \to 0$:
\[
N_\varepsilon (\lambda) = \varepsilon^{-q} \frac{\omega_{n-r}}{{\rm
vol}(V/\Gamma)}
 \sum_{\gamma^*\in\Gamma^*} \left(\frac{\lambda}{4\pi^2}-|\gamma^*-A|^2\right)^{(n-r)/2}
+O(\varepsilon^{\frac{2q}{n-r+1}-q}),
\]
where $\omega_{n-r}$ is the volume of the unit ball in ${\mathbb
R}^{n-r}$.
\end{thm}

\section{Proof of the main results}
\subsection{Proof of Theorem~\ref{t:lattice_points}} We will follow
the proof of Theorem 1.2 in \cite{lattice-points}. Therefore, we
will skip some details, referring the interested reader to
\cite{lattice-points}. First of all, we observe that we have the
inclusion
\[
{\mathbb Z}^n \subset \bigsqcup_{\gamma^*\in \Gamma^*} P_{\gamma^*}.
\]

For any $\gamma^*\in \Gamma^*$, denote
\[
{\mathbb Z}^n_{\gamma^*}={\mathbb Z}^n\bigcap P_{\gamma^*}=\{ k\in
{\mathbb Z}^n : \pi_V(k)=\gamma^*\}.
\]
We identify the affine subspace $P_{\gamma^*}$ with the linear space
$V^\bot$, fixing an arbitrary point $k_{\gamma^*}\in {\mathbb
Z}^n_{\gamma^*}$:
\[
P_{\gamma^*}=k_{\gamma^*}+V^\bot.
\]
It is easy to see that
\[
{\mathbb Z}^n_{\gamma^*}=k_{\gamma^*}+\Gamma^\bot,
\]
where
\[
\Gamma^\bot={\mathbb Z}^n\bigcap V^\bot
\]
is a lattice in $V^\bot$. Observe that
\begin{equation}\label{e:vol}
{\rm vol}(V^\bot/\Gamma^\bot)={\rm vol}(V/\Gamma).
\end{equation}

Thus, we can write
\begin{equation}\label{e:ne}
n_\varepsilon (S)=\sum_{\gamma^*\in \Gamma^*}n_\varepsilon
(S,{\gamma^*}),
\end{equation}
where
\[
n_\varepsilon (S,{\gamma^*})=\# (T_\varepsilon(S)\cap {\mathbb
Z}^n_{\gamma^*}).
\]
Note that, since $S$ is bounded, the sum in the right hand side
of~\eqref{e:ne} has finitely many non-vanishing terms.

Fix $\gamma^*\in \Gamma^*$. Let $\chi_{S_{\gamma^*}}$ be the
indicator of the set $S_{\gamma^*}=S\bigcap P_{\gamma^*}$. It is
easy to see that
\[
n_\varepsilon (S,{\gamma^*}) = \sum_{\gamma\in \Gamma^\bot}
\chi_{S_{\gamma^*}}
(k_{\gamma^*}+(T_{\varepsilon^{-1}}(k_{\gamma^*})-k_{\gamma^*})+T_{\varepsilon^{-1}}(\gamma)).
\]

The space $V^\bot$ decomposes into the direct sum
\begin{equation}\label{e:Vbotdecomp}
V^\bot=F_V\bigoplus H,
\end{equation}
where $F_V=F\cap V^\bot$. We will write the decomposition of $x\in
V^\bot$, corresponding to \eqref{e:Vbotdecomp}, as follows:
\[
x=x_F+x_H, \quad x_F\in F_V, x_H\in H.
\]
Note that
\[
T_{\varepsilon}(x)=x_F +\varepsilon^{-1}x_H.
\]

Let $\rho \in C^\infty_c({\mathbb R})$ be an even function such that
$0\leq \rho(x) \leq 1$ for any $x\in {\mathbb R}$ and ${\rm
supp}\,\rho\subset (-1,1)$. For any $t_F>0$ and $t_H>0$, define a
function $\rho_{t_F,t_H} \in C^\infty_0(V^\bot)$ by
\begin{equation}\label{e:rhotFtH}
\rho_{t_F,t_H}(x)=\frac{c}{t_F^{p-r}t_H^{q}}\rho
\left(\left(t^{-2}_Fx^2_F+t^{-2}_Hx^2_H\right)^{1/2}\right), \quad
x\in V^\bot,
\end{equation}
where the constant $c>0$ is chosen so that $\int_{V^\bot}
\rho_{1,1}(x)\,dx=1$. The function $\rho_{t_F,t_H}$ is supported in
the ellipsoid
\[
B(0,t_F,t_H)=\left\{x\in V^\bot :
\frac{x^2_F}{t^2_F}+\frac{x^2_H}{t^2_H}<1\right\}.
\]

Define the function $n_{\varepsilon,t_F,t_H}(S,\gamma^*)$ by
\[
n_{\varepsilon,t_F,t_H}(S,\gamma^*)=\sum_{k\in {\mathbb
Z}^n_{\gamma^*}}(\chi_{T_\varepsilon(S_{\gamma^*})}\ast
\rho_{t_F,t_H})(k),
\]
where the function $\chi_{T_\varepsilon(S_{\gamma^*})} \ast
\rho_{t_F,t_H}\in C^\infty_0(P_{\gamma^*})$ is defined by
\[
(\chi_{T_\varepsilon(S_{\gamma^*})} \ast \rho_{t_F,t_H})(y)
=\int_{V^\bot} \chi_{T_\varepsilon(S_{\gamma^*})} (y-x)
\rho_{t_F,t_H}(x)\, dx, \quad y\in P_{\gamma^*}.
\]

For any domain $D\subset P_{\gamma^*}$ and for any $t_F>0$ and
$t_H>0$, denote
\[
D_{t_F,t_H}=\bigcup_{x\in D}(x+B(0,t_F,t_H)),
\]
and
\[
D_{-t_F,-t_H} = P_{\gamma^*} \setminus (P_{\gamma^*} \setminus
D)_{t_F,t_H}.
\]
It is easy to see that, for any $\varepsilon>0$, $t_F>0$ and
$t_H>0$, one has
\[
T_\varepsilon(D_{t_F,\varepsilon t_H})=(T_\varepsilon(D))_{t_F,t_H}.
\]

\begin{lem}
For any $\varepsilon>0$, $t_F>0$ and $t_H>0$, we hav
\[
n_{\varepsilon,t_F,t_H}((S_{\gamma^*})_{-t_F,-\varepsilon
t_H},\gamma^*) \leq n_\varepsilon (S,{\gamma^*}) \leq
n_{\varepsilon,t_F,t_H}((S_{\gamma^*})_{t_F,\varepsilon
t_H},\gamma^*).
\]
\end{lem}

For any $f\in {\mathcal S}(V^\bot)$, define its Fourier transform
$\hat{f}\in {\mathcal S}(V^\bot)$ by
\[
\hat{f}(\xi)=\int_{V^\bot}e^{-2\pi i(\xi,x)}f(x)\,dx.
\]

Recall the Poisson summation formula
\begin{equation}\label{e:Poisson}
\sum_{k\in \Gamma^\bot} f(k)=\frac{1}{{\rm
vol}(V/\Gamma)}\sum_{k^*\in {\Gamma^\bot}^*}\hat{f}(k^*), \quad f\in
{\mathcal S}(V^\bot),
\end{equation}
where ${\Gamma^\bot}^*\subset V^{\bot}$ is the dual lattice of
$\Gamma^\bot $, and we used the relation \eqref{e:vol}.

For any $N>0$, we have the estimate
\begin{equation}\label{e:Frho}
|\hat{\rho}_{t_F,t_H}(\xi)|\leq
C_N\frac{1}{1+t_F^N|\xi_F|^N+t_H^N|\xi_H|^N}, \quad \xi \in V^\bot.
\end{equation}
Therefore, we can apply \eqref{e:Poisson} to the function
\begin{equation}\label{e:f}
f(x)=(\chi_{T_\varepsilon((S_{\gamma^*})_{t_F,\varepsilon t_H})}
\ast \rho_{t_F,t_H})(k_{\gamma^*}+x), \quad x\in V^\bot.
\end{equation}
Using the relations
\[
\hat\chi_{T_\varepsilon((S_{\gamma^*})_{t_F,\varepsilon t_H})}(\xi)=
\varepsilon^{-q}e^{2\pi i(\xi,(1-T_{\varepsilon})(k_{\gamma^*}))}
\hat\chi_{(S_{\gamma^*})_{t_F,\varepsilon
t_H}}(T_{\varepsilon}(\xi)).
\]
and
\[
\hat\rho_{t_F,t_H}(\xi)=\hat\rho_{1,1}(t_F\xi_F+t_H\xi_H), \quad
\xi\in V^\bot.
\]
we obtain
\begin{multline}\label{e:ne1}
n_{\varepsilon,t_F,t_H}((S_{\gamma^*})_{t_F,\varepsilon
t_H},\gamma^*)=\sum_{k\in {\mathbb
Z}^n_{\gamma^*}}(\chi_{T_\varepsilon((S_{\gamma^*})_{t_F,\varepsilon
t_H})}\ast \rho_{t_F,t_H})(k) \\ = \frac{\varepsilon^{-q}}{{\rm
vol}(V/\Gamma)} \sum_{k\in {\Gamma^\bot}^*} e^{2\pi
i(k,(1-T_{\varepsilon})(k_{\gamma^*}))}
\hat\chi_{(S_{\gamma^*})_{t_F,\varepsilon t_H}}(
T_{\varepsilon}(k))\hat\rho_{1,1}(t_Fk_F+t_Hk_H).
\end{multline}
We can write
\[
n_{\varepsilon,t_F,t_H}((S_{\gamma^*})_{t_F,\varepsilon
t_H},\gamma^*)=n^{\prime}_{\varepsilon,t_F,t_H}
((S_{\gamma^*})_{t_F,\varepsilon
t_H},\gamma^*)+n^{\prime\prime}_{\varepsilon,t_F,t_H}
((S_{\gamma^*})_{t_F,\varepsilon t_H},\gamma^*),
\]
where
\begin{multline*}
n^{\prime}_{\varepsilon,t_F,t_H} ((S_{\gamma^*})_{t_F,\varepsilon
t_H},\gamma^*)\\ = \frac{\varepsilon^{-q}}{{\rm vol}(V/\Gamma)}
\sum_{k\in {\Gamma^\bot}^*, k_H=
0}\hat\chi_{(S_{\gamma^*})_{t_F,\varepsilon
t_H}}(T_{\varepsilon}(k))\hat\rho_{1,1}(t_Fk_F+t_Hk_H),
\end{multline*}
and
\begin{multline*}
n^{\prime\prime}_{\varepsilon,t_F,t_H}
((S_{\gamma^*})_{t_F,\varepsilon t_H},\gamma^*)\\ =
\frac{\varepsilon^{-q}}{{\rm vol}(V/\Gamma)} \sum_{k\in
{\Gamma^\bot}^*, k_H\neq 0} e^{2\pi i(1-\varepsilon^{-1})
(k_H,k_{\gamma^*})}\hat\chi_{(S_{\gamma^*})_{t_F,\varepsilon t_H}}(
T_{\varepsilon}(k)) \hat\rho_{1,1}(t_Fk_F+t_Hk_H).
\end{multline*}

Let $k\in {\Gamma^\bot}^*$ be such that $k_H=0$. Then $k\in F_V$.
Since ${\Gamma^\bot}^*\subset {\mathbb Q}^n$ and $F_V \cap {\mathbb
Q}^n=\{0\}$, we get $k=0$. Thus, we have
\begin{align*}
n^{\prime}_{\varepsilon,t_F,t_H} ((S_{\gamma^*})_{t_F,\varepsilon
t_H},\gamma^*) = & \frac{\varepsilon^{-q}}{{\rm vol}(V/\Gamma)}{\rm
vol}_{n-r} (P_{\gamma^*}\cap S)\\ & +\frac{\varepsilon^{-q}}{{\rm
vol}(V/\Gamma)}{\rm vol}_{n-r} ((S_{\gamma^*})_{t_F,\varepsilon t_H}
\setminus S_{\gamma^*}).
\end{align*}
Since
\[
{\rm vol}_{n-r} ((S_{\gamma^*})_{t_F,\varepsilon t_H} \setminus
S_{\gamma^*})\leq C (t_F+t_H\varepsilon),
\]
we obtain that
\begin{equation}\label{e:n11}
n^{\prime}_{\varepsilon,t_F,t_H}
((S_{\gamma^*})_{t_F,\varepsilon t_H},\gamma^*) =
\frac{\varepsilon^{-q}}{{\rm vol}(V/\Gamma)}{\rm vol}_{n-r}
(P_{\gamma^*}\cap S)+O(t_F\varepsilon^{-q}+t_H\varepsilon^{1-q}).
\end{equation}

Consider the case when $k\in {\Gamma^\bot}^*$ and $k_H\neq 0$. Here
the arguments depend on the conditions on the domain $S$ we have.

(1) Assume that for any $x\in F$ the domain $S\cap \{x+H\}$ is
strictly convex. For any $t\in F_V$ and for any domain $D\subset
P_{\gamma^*}$, we denote
\[
D(t)=\{x_H\in H : k_{\gamma^*}+t+x_H\in D\}\subset H.
\]
For any function $\phi\in {\mathcal S}(H)$, denote by $F_H(\phi)\in
{\mathcal S}(H)$ its Fourier transform:
\[
[F_H(\phi)](\xi_H)=\int_H \phi(x_H) e^{-2\pi i(\xi_H,x_H)}\,dx_H,
\quad \xi_H\in H.
\]
It is easy to see that
\[
\hat\chi_{(S_{\gamma^*})_{t_F,\varepsilon
t_H}}(T_{\varepsilon}(k))=\int_{F_V} e^{-2\pi i(k_F,x_F)}
F_H[\chi_{(S_{\gamma^*})_{t_F,\varepsilon t_H}( x_F)}]
(\varepsilon^{-1}k_H) dx_F,
\]
and, therefore.
\begin{equation}\label{e:chi}
|\hat\chi_{(S_{\gamma^*})_{t_F,\varepsilon
t_H}}(T_{\varepsilon}(k))|\leq \int_{F_V}
|F_H[\chi_{(S_{\gamma^*})_{t_F,\varepsilon t_H}( x_F)}]
(\varepsilon^{-1}k_H)| dx_F.
\end{equation}

By assumption, the domain $S_{\gamma^*}(x_F)=S\cap
\{\gamma^*+x_F+H\}$ is strictly convex. Therefore, for any
sufficiently small $\varepsilon>0$, $t_F>0$ and $t_H>0$, the domain
$(S_{\gamma^*})_{t_F,\varepsilon t_H}(x_F)$ is strictly convex. We
have the estimate
\begin{equation}\label{e:FH2}
|F_H[\chi_{(S_{\gamma^*})_{t_F,\varepsilon t_H}(x_F)}](\xi)| =
O(|\xi|^{-(q+1)/2}), \quad |\xi|\to \infty.
\end{equation}
Thus, using the estimate \eqref{e:chi}, \eqref{e:FH2} and
\eqref{e:Frho}, we obtain that
\begin{multline}\label{e:n22}
|n^{\prime\prime}_{\varepsilon,t_F,t_H} ((S_{\gamma^*})_{t_F,\varepsilon t_H},\gamma^*)|\\
\begin{aligned} & \leq  C \varepsilon^{-q} \sum_{k\in {\Gamma^\bot}^*, k_H\neq
0}\varepsilon^{(q+1)/2} |k_H|^{-(q+1)/2} \frac{1}{1+t_F^N|k_F|^N+t_H^N|k_H|^N}\\
& \leq  C \varepsilon^{-q} \varepsilon^{(q+1)/2}\int_{V^{\bot}}
|x_H|^{-(q+1)/2}\frac{dx_F\,dx_H}{1+t_F^N|x_F|^N+t_H^N|x_H|^N}\\
& \leq C \varepsilon^{-(q-1)/2} t_F^{-(p-r)}t_H^{-(q-1)/2}.
\end{aligned}
\end{multline}
Put $t_F=\varepsilon^{\alpha_F}$, $t_H=\varepsilon^{\alpha_H}$,
where
\[
\alpha_F=\frac{2q}{q+1+2(p-r)}, \quad
\alpha_H=\frac{q-1-2(p-r)}{q+1+2(p-r)}.
\]
Using the estimates \eqref{e:n11} and \eqref{e:n22}, we immediately
conclude the proof of the statement (1).

(2) Assume that $P_{\gamma^*} \cap S$ is strictly convex. Then we
have
\[
|\hat\chi_{(S_{\gamma^*})_{t_F,\varepsilon
t_H}}(T_{\varepsilon}(k))| \leq
C(|k_F|+|\varepsilon^{-1}k_H|)^{-(n-r+1)/2}.
\]
Using this fact, as in \eqref{e:n22}, we obtain
\begin{equation}\label{e:n22a}
|n^{\prime\prime}_{\varepsilon,t_F,t_H}
((S_{\gamma^*})_{t_F,\varepsilon t_H},\gamma^*)| \leq C
\varepsilon^{-q} \varepsilon^{(n-r+1)/2} t_F^{-(p-r)} t_H^{-q}
t_H^{(n-r+1)/2}.
\end{equation}
To complete the proof of the statement (2), we put
$t_F=\varepsilon^{\alpha_F}$, $t_H=\varepsilon^{\alpha_H}$, where
$\alpha_F\geq 0$ and $\alpha_H\geq -1$ are given by
\[
\alpha_F=\frac{2q}{n-r+1}, \quad \alpha_H=\frac{q-p+r-1}{n-r+1}.
\]

\subsection{Proof of Theorem~\ref{t:VarchenkoH}} Let
$h\in {\rm SO}(H)$.  We apply the formula~\eqref{e:ne1} to the set
$h(S)$. Since $T_\varepsilon h=hT_\varepsilon $, we have
\[
(h(S)_{\gamma^*})_{t_F,\varepsilon
t_H}=(h(S_{\gamma^*}))_{t_F,\varepsilon
t_H}=h((S_{\gamma^*})_{t_F,\varepsilon t_H}).
\]
Therefore, the formula \eqref{e:ne1} reads as
\begin{multline} \label{e:ne1bis}
n_{\varepsilon,t_F,t_H}((h(S)_{\gamma^*})_{t_F,\varepsilon
t_H},\gamma^*)\\ =n^{\prime}_{\varepsilon,t_F,t_H}
((h(S)_{\gamma^*})_{t_F,\varepsilon t_H},\gamma^*)
+n^{\prime\prime}_{\varepsilon,t_F,t_H}
((h(S)_{\gamma^*})_{t_F,\varepsilon t_H},\gamma^*),
\end{multline}
where, since $h$ preserves the volume in ${\mathbb R}^n$, the first
term is independent of $h$ and satisfies the estimate
\[
n^{\prime}_{\varepsilon,t_F,t_H} ((h(S)_{\gamma^*})_{t_F,\varepsilon
t_H},\gamma^*) =\frac{\varepsilon^{-q}}{{\rm vol}(V/\Gamma)}{\rm
vol}_{n-r} (P_{\gamma^*}\cap
S)+O(t_F\varepsilon^{-q}+t_H\varepsilon^{1-q}),
\]
and
\begin{multline*}
n^{\prime\prime}_{\varepsilon,t_F,t_H}
((h(S)_{\gamma^*})_{t_F,\varepsilon t_H},\gamma^*) =
\frac{\varepsilon^{-q}}{{\rm vol}(V/\Gamma)} \sum_{k\in
{\Gamma^\bot}^*, k_H\neq 0} e^{2\pi i
(k_H,(1-\varepsilon^{-1}h)k_{\gamma^*})} \times \\
\times \hat\chi_{(S_{\gamma^*})_{t_F,\varepsilon
t_H}}(h^tT_{\varepsilon}(k)) \hat\rho_{1,1}(t_Fk_F+t_Hk_H).
\end{multline*}
Here $h^t=h^{-1}$ denotes the transpose of $h$.

Consider the case when $k\in {\Gamma^\bot}^*$ and $k_H\neq 0$. We
will keep notation used in the previous subsection. As in
\eqref{e:chi}, we have
\[
|\hat\chi_{(S_{\gamma^*})_{t_F,\varepsilon
t_H}}(h^tT_{\varepsilon}(k))|\leq \int_{F_V}
|F_H[\chi_{[(S_{\gamma^*})_{t_F,\varepsilon t_H}]( x_F)}]
(\varepsilon^{-1}h^tk_H)| dx_F.
\]
Using the results of \cite{Brandolini}, we get
\begin{multline*}
\int_{{\rm SO}(H)}|F_H[\chi_{[(S_{\gamma^*})_{t_F,\varepsilon t_H}](
x_F)}] (\varepsilon^{-1}h^tk_H)|\,dh\\
\leq \left(\int_{{\rm
SO}(H)}|F_H[\chi_{[(S_{\gamma^*})_{t_F,\varepsilon t_H}]( x_F)}]
(\varepsilon^{-1}h^tk_H)|^2\,dh\right)^{1/2}
=O(\varepsilon^{(q+1)/2}).
\end{multline*}
Proceeding as above (cf. \eqref{e:n22}), we immediately complete the
proof.

\subsection{Proof of Theorem~\ref{t:Varchenko-full}} Now we assume
that $h\in {\rm SO}(V^\bot)$. Then we still have
\[
h(S)_{\gamma^*}=h(S_{\gamma^*}),
\]
but the equality
\[
(h(S_{\gamma^*}))_{t_F,\varepsilon
t_H}=h((S_{\gamma^*})_{t_F,\varepsilon t_H})
\]
only if
\[
t_F=\varepsilon t_H=t.
\]
Therefore, the formula \eqref{e:ne1bis} holds for such $t_F$ and
$t_H$.

Let $k\in {\Gamma^\bot}^*$ with $k_H\neq 0$. Then we have
\begin{align*}
\hat\chi_{(S_{\gamma^*})_{t_F,\varepsilon
t_H}}(h^tT_{\varepsilon}(k)) = &
\int_{V^\bot}\chi_{(S_{\gamma^*})_{t_F,\varepsilon
t_H}}(k_{\gamma^*}+x)e^{-2\pi i(h^tT_{\varepsilon}(k),x)}dx\\
= & \int_{D_\varepsilon} e^{-2\pi
i(T_{\varepsilon}(k),h(x-k_{\gamma^*}))}\,dx.
\end{align*}
where, for simplicity of notation, we put
$D_\varepsilon=(S_{\gamma^*})_{t_F,\varepsilon t_H}$. By Stokes'
formula, we obtain
\[
\hat\chi_{(S_{\gamma^*})_{t_F,\varepsilon
t_H}}(h^tT_{\varepsilon}(k)) =  - \frac{1}{2\pi
i|T_{\varepsilon}(k)|}\int_{\partial D_\varepsilon} e^{-2\pi
i(T_{\varepsilon}(k),h(x-k_{\gamma^*}))}
i_{h^t\mathbf{n}_\varepsilon}dx,
\]
where $\mathbf{n}_\varepsilon
=\frac{T_{\varepsilon}(k)}{|T_{\varepsilon}(k)|}$. We can write
\begin{multline*}
\int_{{\rm SO}(V^\bot)} |\hat\chi_{(S_{\gamma^*})_{t_F,\varepsilon
t_H}}(h^tT_{\varepsilon}(k))|^2\,dh\\
= \frac{1}{4\pi^2|T_{\varepsilon}(k)|^2}\int_{\partial
D_\varepsilon}\int_{\partial D_\varepsilon}\int_{{\rm SO}(V^\bot)}
e^{-2\pi i(k_F+\varepsilon^{-1}k_H,h(x-y))}
i_{h^t\mathbf{n}_\varepsilon}dx\wedge
i_{h^t\mathbf{n}_\varepsilon}dy \wedge dh.
\end{multline*}

So this is an oscillating integral with the phase
\[
\Phi(x,y,h)=(k_H,h(x-y)), \quad x,y\in \partial D_\varepsilon
\subset V^\bot, h\in {\rm SO}(V^\bot).
\]
If $(x_0,y_0,h_0)\in \partial D_\varepsilon \times \partial
D_\varepsilon\times {\rm SO}(V^\bot)$ is a critical point of $\Phi$,
then we have:
\begin{itemize}
  \item for any $v\in T_{x_0}\partial D_\varepsilon$ and $w\in T_{y_0}\partial D_\varepsilon$
  \[
(k_H,h_0v)=(h^t_0k_H,v)=0, \quad (k_H,h_0w)=(h^t_0k_H,w)=0;
  \]
  \item for any $X\in \mathfrak{so}(V^\bot)$
\begin{equation}\label{e:crit2}
(k_H,h_0X(x_0-y_0))=(h^t_0k_H,X(x_0-y_0))=0.
  \end{equation}
\end{itemize}
By \eqref{e:crit2}, it follows that $x_0-y_0=\alpha h^t_0k_H$ with
some $\alpha\in \mathbb R$.

As in \cite[Lemma 1]{Varchenko83}, we have that at any critical
point of $\Phi$ on $\partial D_\varepsilon\times \partial
D_\varepsilon\times {\rm SO}(V^\bot)$ the rank of its second
differential is at least $2(n-r)-2$, that implies by a slight
modification of \cite[Lemma 2]{Varchenko83} that
\[
\int_{{\rm SO}(V^\bot)} |\hat\chi_{(S_{\gamma^*})_{t_F,\varepsilon
t_H}}(h^tT_{\varepsilon}(k))|^2\,dh\leq
C|\varepsilon^{-1}k_H|^{-(n-r+1)},
\]
and therefore
\[
\int_{{\rm SO}(V^\bot)} |\hat\chi_{(S_{\gamma^*})_{t_F,\varepsilon
t_H}}(h^tT_{\varepsilon}(k))|\,dh\leq
C|\varepsilon^{-1}k_H|^{-(n-r+1)/2},
\]
As in \eqref{e:n22a}, we obtain
\begin{multline*}
\int_{{\rm SO}(V^\bot)} |n^{\prime\prime}_{\varepsilon,t_F,t_H}
((h(S)_{\gamma^*})_{t_F,\varepsilon t_H},\gamma^*)|\,dh \\
\leq C
\varepsilon^{-q} \varepsilon^{(n-r+1)/2} t_F^{-(p-r)} t_H^{-q}
t_H^{(n-r+1)/2} = Ct^{-(n-r-1)/2}.
\end{multline*}
To complete the proof, we put $t=\varepsilon^{\alpha}$, where
\[
\alpha=\frac{2q}{n-r+1}.
\]

\subsection{Proof of Theorem~\ref{t:Varchenko-full2}} First,
we observe that there exists an invariant measure $dh$ on the
homogeneous space $SO(V^\bot)\backslash SO(n)$ such that
\[
\int_{SO(n)} |R_\varepsilon (hS)|\,dh = \int_{SO(V^\bot)\backslash
SO(n)}\left( \int_{SO(V^\bot)} |R_\varepsilon
(h_1h_2S)|\,dh_1\right)dh_2.
\]
For any $h_2\in SO(n)$, by Theorem~\ref{t:Varchenko-full} applied to
the domain $h_2S$, we have
\[
\int_{SO(V^\bot)} |R_\varepsilon (h_1h_2S)|\,dh_1\leq
C\varepsilon^{\frac{2q}{n-r+1}-q},\quad \varepsilon > 0.
\]
Moreover, it is easy to see from \cite[Lemma 2]{Varchenko83} that
the constant $C>0$ can be chosen independent of $h_2\in SO(n)$, that
immediately completes the proof of Theorem~\ref{t:Varchenko-full2}.

\section*{Acknowledgements} This research was partly supported by
the Russian Foundation of Basic Research (grants no. 09-01-00389 and
no. 12-01-00519).

\end{document}